\documentclass[leqno,11pt]{amsart}

\usepackage{amsfonts}
\usepackage{amsmath}
\usepackage{amsfonts}
\usepackage{amssymb}
\usepackage{amscd}

\newtheorem{theorem}{Theorem}
\newtheorem*{akn}{Acknowledgements}

\newtheorem{corollary}[theorem]{Corollary}

\newtheorem*{example*}{Example}

\newtheorem{lemma}[theorem]{Lemma}

\newtheorem{prop}[theorem]{Proposition}
\newtheorem{remark}[theorem]{Remark}

\newcommand{\erre}{\mathbb{R}}

\newcommand{\enne}{\mathbb{N}}

\newcommand{\ricc}{\operatorname{Ric}}

\newcommand{\hess}{\operatorname{hess}}
\newcommand{\Hess}{\operatorname{Hess}}

\newcommand{\R}{\operatorname{R}}

\newcommand{\ra}{\rightarrow}

\newcommand{\norm}[1]{{\left\|#1\right\|}}              
\newcommand{\pair}[1]{\left\langle#1\right\rangle}      

\renewcommand{\tilde}[1]{\widetilde{#1}}

\begin{document}
\title[Comparison Theorems in Lorentzian Geometry]{Comparison Theorems in Lorentzian Geometry and applications to spacelike hypersurfaces}
\author[Debora Impera]{Debora Impera}
\address{Dipartimento di Matematica,
Universit\`a
degli studi di Milano, via Saldini 50, I-20133 Milano, Italy.}
\email{debora.impera@unimi.it}
\date{\today}

\begin{abstract}
In this paper we prove Hessian and Laplacian comparison theorems for the Lorentzian distance function in a spacetime with sectional (or Ricci) curvature bounded by a certain function by means of a comparison criterion for Riccati equations. Using these results, under suitable conditions, we are able to obtain some estimates on the higher order mean curvatures of spacelike hypersurfaces satisfying a Omori-Yau maximum principle for certain elliptic operators.
\end{abstract}

\maketitle

\section{Introduction}
In general relativity each point of a Lorentzian manifold corresponds to an event. The events that we may experience in the universe are the ones in our chronological future, hence it may be interesting to investigate the geometry of this one. This can be done by means of the analysis of the Lorentzian distance function. Unfortunately this function is not differentiable in any spacetime; precisely, it is not even continuous in general. Nevertheless, in strongly causal spacetimes, the Lorentzian distance function from a point is differentiable at least in a ``sufficiently near'' chronological future of each point. In this case is possible to analyze the geometry of spacetimes by means of the level sets of the Lorentzian distance function with respect to this point. To do that, the main tools are Hessian and Laplacian comparison theorems for the Lorentzian distance of the spacetime, hence many works have been written in this spirit.\\ 
For instance, in a recent paper by F. Erkekoglu, E. Garc\'ia-Rio and D. N. Kupeli (\cite{egrk}), following the approach of R. E. Greene and H. Wu in \cite{greenewu}, the authors obtain Hessian and Laplacian comparison theorems for the Lorentzian distance functions of Lorentzian manifolds comparing their sectional curvatures.
Afterwards, in \cite{aliashurtadopalmer}, L. J. Al\'ias, A. Hurtado and V. Palmer use these theorems to study the Lorentzian distance function restricted to a spacelike hypersurface $\Sigma^n$ immersed into a spacetime $M^{n+1}$. In particular, under suitable conditions, they derive sharp estimates for the mean curvature of spacelike hypersurfaces with bounded image in the ambient spacetime.\\
In this paper we obtain Hessian and Laplacian comparison theorems for Lorentzian manifolds with sectional curvature of timelike planes bounded by a function of the Lorentzian distance, improving in this way on classical results, and we give some applications to the study of spacelike hypersurfaces. 

The paper is organized as follows. In Sections 2 and 3 we present some basic concepts and terminology involving the Lorentzian distance function from a point and we prove our Hessian and Laplacian comparison theorems. To obtain these theorems we use an `analytic' approach inspired by P. Petersen (\cite{petersen}) avoiding, in this way, the `geometric' approach used by Greene and Wu. 
In Section 4 we focus on the study of the Lorentzian distance function restricted to spacelike hypersurfaces. Hence, using the Omori-Yau maximum principle, we derive some estimates on the mean curvature that generalize the ones in \cite{aliashurtadopalmer}. Moreover, using a generalized Omori-Yau maximum principle for certain elliptic operators, we also obtain some estimates for the higher order mean curvatures associated to the immersion. Finally, in Section 5, we restrict ourselves to the case when the ambient space has constant sectional curvature and we prove a Bernstein-type theorem for spacelike hypersurfaces with constant $k$-mean curvature that generalizes Corollary 4.6 in \cite{aliashurtadopalmer}.  
  
\section{Preliminaries}

Let $M^{n+1}$ be an $n+1$-dimensional spacetime, that is, an $n+1$-dimensional time-oriented Lorentzian manifold and let $p,\ q\in M$. Using the standard terminology and notation in Lorentzian geometry, we say that $q$ is in the chronological future of $p$, written $p\ll q$, if there exists a future-directed timelike curve from $p$ to $q$. Similarly, we say that $q$ is in the causal future of $p$, written $p<q$, if there exists a future-directed causal (that is nonspacelike) curve from $p$ to $q$. 
For a subset $S\subset M$, we define the chronological future of $S$ as
\begin{displaymath}
I^+(S)=\lbrace q\in M |p\ll q \text{ for some }p\in S \rbrace,
\end{displaymath}
and  the causal future of $S$ as
\begin{displaymath}
J^+(S)=\lbrace q\in M |p\leq q \text{ for some }p\in S \rbrace,
\end{displaymath}
where $p\leq q$ means that either $p<q$ or $p=q$.
In particular, the chronological and the causal future of a point $p\in M$ are, respectively
\begin{displaymath}
I^+(p)=\lbrace q\in M |p\ll q \rbrace,\qquad J^+(p)=\lbrace q\in M |p\leq q \rbrace.
\end{displaymath}
 It is well known that $I^+(p)$ is always open, while $J^+(p)$ is neither open nor closed in general.
Let $q \in J^+(p)$. Then the Lorentzian distance $d(p,q)$ is defined as the supremum of the Lorentzian lengths of all the future-directed causal curves from $p$ to $q$. If $q\not\in J^+(p)$, then $d(p,q)=0$ by definition. Moreover, $d(p,q)>0$ if and only if $q\in J^+(p)$.
 Given a point $p\in M$ one can define the Lorentzian distance function $d_p: M\ra[0,+\infty)$ with respect to 
$p$ by
$$
d_p(q)=d(p,q).
$$
Let
$$
T_{-1}M_{|p}=\lbrace v\in T_pM|v \text{ is a future-directed timelike unit vector}\rbrace
$$
be the fiber of the unit future observer bundle of $M^{n+1}$ at $p$. Define the function
$$
s_p:T_{-1}M_{|p}\ra [0,+\infty],\qquad s_p(v)=\sup \lbrace t \geq 0\ |\ d_p(\gamma_v(t))=t\rbrace,
$$
where $\gamma_v:[0,a)\ra M$ is the future timelike geodesic with $\gamma_v(0)=p$, $\gamma_v'(0)=v$. The future timelike cutlocus $\Gamma^+(p)$ of $p$ in $T_pM$ is defined as
$$
\Gamma^+(p)=\{s_p(v)v\ |\ v\in T_pM \text{ and }0<s_p(v)<+\infty\}
$$
and the future timelike cutlocus $C^+_t(p)$ of $p$ in $M$ is $C^+_t(p)=\exp_p(\Gamma^+(p))$ wherever the exponential map $\exp_p$ at $p$ is defined on $\Gamma^+(p)$.\\
It is well known that the Lorentzian distance function on arbitrary spacetimes may fail in general to be continuous and finite valued. It is known that this is true for globally hyperbolic spacetimes. We recall that a spacetime $M$ is said to be \textit{globally hyperbolic} if it is strongly causal and it satisfies the condition that $J^+(p)\cap J^-(q)$ is compact for all $p,\ q\in M$. Moreover, a Lorentzian manifold $M$ is said to be \textit{strongly causal} at a point $p\in M$ if for any neighborhood $U$ of $p$ there exists no timelike curve that passes through $U$ more than once. In general, in order to guarantee the smoothness of this function we need to restrict it on certain special subsets of $M$. 
Let  
$$
\tilde{\mathcal{I}}^+(p)=\{tv\ | \ v\in T_{-1}M_{|p} \text{ and } 0<t<s_{p}(v)\},
$$
and let
$$
\mathcal{I}^+(p)=\mathrm{exp}(\mathrm{int}(\tilde{\mathcal{I}}^+(p)))\subset I^+(p).
$$
Since
$$
\exp_p:\rm{int}(\tilde{\mathcal{I}}^+(p))\ra\mathcal{I}^+(p)
$$
is a diffeomorphism, $\mathcal{I}^+(p)$ is an open subset of $M$.
In the lemma below we summarize the main properties of the Lorentzian distance function. \begin{lemma}[\cite{egrk}, Section 3.1]
Let $M$ be a spacetime and $p\in M$. 
\begin{enumerate}
\item If $M$ is strongly causal at $p$, then $s_p(v)>0\ \forall v\in T_{-1}M_{|p}$ and $\mathcal{I}^+(p) \neq \emptyset$,\\
\item If $\mathcal{I}^+(p)\neq \emptyset$, then the Lorentzian distance function $d_p$ is smooth on $\mathcal{I}^+(p)$ and $\overline{\nabla}d_p$ is a past-directed timelike (geodesic) unit vector field on $\mathcal{I}^+(p)$.
\end{enumerate}
\end{lemma}
\begin{remark}
\rm{
If $M$ is a globally hyperbolic spacetime and $\Gamma^+(p)=\emptyset$, then $\mathcal{I}^+(p)=I^+(p)$ and hence the Lorentzian distance function $d_p$ with respect to $p$ is smooth on $I^+(p)$ for each $p\in M$.\\
We also observe that if $M$ is a Lorentzian space form, then it is globally hyperbolic and geodesically complete. Moreover, every timelike geodesic realizes the distance between its points. Hence $\Gamma^+(p)=\emptyset$ and we conclude again that the Lorentzian distance function $d_p$ is smooth on $I^+(p)$ for each $p\in M$. 
}
\end{remark}
\section{Hessian and Laplacian Comparison Theorems}
This section is devoted to exhibit estimates for the Hessian and the Laplacian of the Lorentzian distance function in Lorentzian manifolds under conditions on the sectional or Ricci curvature.
To prove our theorems we will need the following Sturm comparison result.
\begin{lemma}\label{sturm}
Let $G$ be a continuous function on $[0,+\infty)$ and let $\phi,\ \psi \in C^1([0,+\infty))$ with $\phi',\ \psi' \in AC([0,+\infty))$ be solutions of the problems
\begin{displaymath}
\left\{\begin{array}{l}
\phi''-G\phi\leq 0 \qquad \text{a.e. in } (0,+\infty)\\
\phi(0)=0
\end{array}\right.
\left\{\begin{array}{l}
\psi''-G\psi\geq 0 \qquad \text{a.e. in } (0,+\infty)\\
\psi(0)=0,\ \psi'(0)>0
\end{array}\right.
\end{displaymath}
If $\phi(r)>0$ for $r\in(0,T)$ and $\psi'(0)\geq \phi'(0)$, then $\psi(r)>0$ in $(0,T)$ and
$$
\frac{\phi'}{\phi}\leq \frac{\psi'}{\psi} \text{ and }\psi \geq \phi \qquad \text{on } (0,T). 
$$
\end{lemma}
For a proof of the lemma see \cite{pirise2}.
Using the Sturm comparison result, we obtain a comparison result for solutions of Riccati inequalities with appropriate asymptotic behaviour. 
\begin{corollary}\label{riccati}
Let $G$ be a continuous function on $[0,+\infty)$ and let $g_i\in AC((0,T_i))$ be solutions of the Riccati differentials inequalities
\begin{displaymath}
g_1'-\frac{g_1^2}{\alpha}+\alpha G\geq 0,\ (\mathrm{ resp. }\leq 0)\qquad g_2'+\frac{g_2^2}{\alpha}-\alpha G\geq 0,\ (\mathrm{ resp. }\leq 0)
\end{displaymath}
a.e. in $(0,T_i)$, satisfying the asymptotic conditions
\begin{displaymath}
g_i(t)=\frac{\alpha}{t}+o(t)\qquad \text{as}\qquad t\ra 0^+, \end{displaymath}
for some $\alpha>0$. Then $T_1 \leq T_2$ (resp. $T_1 \geq T_2$) and $-g_1(t)\leq g_2(t)$ in $(0,T_1)$ (resp. $-g_2(t)\leq g_1(t)$ in $(0,T_2)$).  
\end{corollary}
\begin{proof}
Since $\tilde{g_i}=\alpha^{-1}g_i$ satisfies the conditions in the statement with $\alpha=1$, without loss on generality we may assume 
that $\alpha=1$. Notice that $g_i(s)-\frac{1}{s}$ is bounded and integrable in a neighbourhood of $s=0$. Hence the same is true for 
the function $-g_1(s)-\frac{1}{s}$. Indeed
$$
-\Big(g_1(s)+\frac{1}{s}\Big)<-\Big(g_1(s)-\frac{1}{s}\Big)\leq \Big|g_1(s)-\frac{1}{s}\Big|\leq C,
$$
for some constant $C>0$.
Now let $\phi_i\in C^1([0,T_i))$ be the positive functions defined by
$$
\phi_1(t)=t\exp\Big(-\int_0^t\big(g_1(s)+\frac1s\big)ds\Big),\ \phi_2(t)=t\exp\Big(\int_0^t\big(g_2(s)-\frac1s\big)ds\Big).
$$  
Then $\phi_i(0)=0$, $\phi_i'\in AC((0,T_i))$, $\phi_i'(0)=1$ and
$$
\phi_1'(t)=-g_1(t)\phi_1(t),\ \phi_2'(t)=g_2(t)\phi_2(t)
$$
Hence
$$
\phi_1''\leq G\phi_1,\quad \phi_2''\geq G\phi_2 \quad(\mathrm{resp. }\  \phi_1''\geq G\phi_1,\quad \phi_2''\leq G\phi_2).
$$
Then, it follows by Lemma \ref{sturm} that $T_1 \leq T_2$ (resp. $T_1 \geq T_2$) and
$$
-g_1(t)=\frac{\phi_1'(t)}{\phi_1(t)}\leq \frac{\phi_2'(t)}{\phi_2(t)}=g_2(t)
\quad(\mathrm{resp. } -g_2(t)=\frac{\phi_2'(t)}{\phi_2(t)}\leq \frac{\phi_1'(t)}{\phi_1(t)}=g_1(t)).
$$
\end{proof}
We are now ready to prove the Hessian and Laplacian comparison theorems. In both cases we will follow the proofs given by S. Pigola, 
M. Rigoli and A. G. Setti in \cite{pirise2} of the corresponding theorems in the Riemannian setting.\\ We will denote by 
$\overline{\nabla}$ and $\overline{\Delta}$ respectively the Levi-Civita connection and the 
Laplacian on the spacetime $M$. Moreover, for a given function
$f\in C^{2}(M)$,
we denote by $\overline{\hess}{f}:TM\rightarrow TM$ the symmetric operator given by 
$\overline{\hess}{f}(X)=\overline{\nabla}_X\overline{\nabla} f$ for every $X\in TM$, and by
$\overline{\Hess}{f}:TM\times TM\rightarrow C^{\infty}(M)$ the metrically equivalent bilinear form given by
\[
\overline{\Hess}{f}(X,Y)=\pair{\overline{\hess}{f}(X),Y}.
\] 
\begin{theorem}[Hessian Comparison Theorem]
Let $M^{n+1}$ be an $n+1$-dimensional spacetime. Assume that there exists a point $p \in M$ such that $\mathcal{I}^+(p)\neq\emptyset$ 
and let $r(\cdot)=d_p(\cdot)$ be the Lorentzian distance function from $p$. Given a smooth even 
function $G$ on $\erre$, let $h$ be a solution of the Cauchy problem
\begin{displaymath}
\left\{\begin{array}{l}
h''-Gh= 0\\
h(0)=0,\ h'(0)=1
\end{array}\right.
\end{displaymath}
and let $I=[0,r_0)\subset[0,+\infty)$ be the maximal interval where $h$ is positive and $q \in \mathcal{I}^+(p)\cap B^+(p,r_0)$, 
where 
\begin{displaymath}
B^+(p,r_0)=\lbrace q\in I^+(p) |d_p(q)<r_0\rbrace.
\end{displaymath}
If
\begin{equation}\label{curvupbound}
K_{M}(\Pi)\leq G(r)
\end{equation}
for all timelike planes $\Pi$, then
\begin{displaymath}
\overline{\Hess} r(X,X) \geq -\frac{h'}{h}(r)\pair{X,X}
\end{displaymath}
for every spacelike $X \in T_qM$ which is orthogonal to $\overline{\nabla}r$.
Analogously, if   
\begin{equation}\label{curvlowbound}
K_{M}(\Pi)\geq G(r)
\end{equation}
for all timelike planes $\Pi$, then
\begin{displaymath}
\overline{\Hess} r(X,X) \leq -\frac{h'}{h}(r)\pair{X,X}
\end{displaymath}
for every spacelike $X \in T_qM$ which is orthogonal to $\overline{\nabla}r$.
\end{theorem}
\begin{proof}
Let $v\in\exp_p^{-1}(q)\in\rm{int}(\tilde{\mathcal{I}}^+(p))$ and let $\gamma(t)=\exp_p(tv)$, $0\leq t\leq s_p(v)$, be the radial 
future directed unit timelike geodesic with $\gamma(0)=p$, $\gamma(s)=q$, $s=r(q)$. Recall that $\gamma'(s)=-\overline{\nabla} r(q)$ 
and $\overline{\nabla}_{\overline{\nabla} r}\overline{\nabla} r(q)=0$. Since $\overline{\nabla} r$ satisfies the timelike eikonal 
inequality, $\overline{\Hess} r$ is diagonalizable (see \cite{grk} Chapter 6 or \cite{grk2} for more details) and $T_qM$ has an 
orthonormal basis consisting of eigenvectors of $\overline{\Hess} r$. Let us denote by $\lambda_{\max}(q)$ and $\lambda_{\min}(q)$ 
respectively its greatest and smallest eigenvalues in the orthogonal complement of $\overline{\nabla} r(q)$. Notice that 
the theorem is proved once one shows that 
\begin{itemize}
\item[$(a)$] if \eqref{curvupbound} holds, then
$$
\lambda_{\min}(q)\geq -\frac{h'}{h}(r(q)).
$$
\item[$(b)$] if \eqref{curvlowbound} holds, then
$$
\lambda_{\max}(q)\leq -\frac{h'}{h}(r(q)).
$$
\end{itemize}   
Let us prove claim $(a)$ first. We claim that if \eqref{curvupbound} holds, then $\lambda_{\min}$ satisfies
\begin{equation}\label{lambdamin}
\left\{\begin{array}{l}
\frac{d}{dt}(\lambda_{\min}\circ\gamma)-(\lambda_{\min}\circ\gamma)^2\geq-G \qquad\text{for a.e. }t>0\\
\lambda_{\min}\circ\gamma=\frac1t+o(t)\qquad\text{as }t\ra0^+ 
\end{array}
\right.
\end{equation}
Namely, by the definition of covariant derivative
\begin{displaymath}
(\overline{\nabla}_X\overline{\hess} u)(Y)=\overline{\nabla}_X(\overline{\hess} u(Y))-\overline{\hess} u (\overline{\nabla}_X Y).
\end{displaymath}
Hence, recalling the definition of the curvature tensor we find  
\begin{displaymath}
(\overline{\nabla}_Y \overline{\hess} u)(X)-(\overline{\nabla}_X \overline{\hess} u)(Y)=\overline{\R}
(X,Y)\overline{\nabla} u.
\end{displaymath}
Choose $u=r$, $X=\overline{\nabla} r$. For every spacelike unit vector $Y\in T_qM$, $Y$ is orthogonal to $\gamma'(s)$ and we can define a vector field $Y$ orthogonal to $\gamma'$ by parallel translation along $\gamma$. Then
\begin{align*}
\overline{\nabla}_{\gamma'(s)}(\overline{\hess} r(Y))=&(\overline{\nabla}_{\gamma'(s)}\overline{\hess} r)(Y)+\overline{\hess} r(\overline{\nabla}_{\gamma'(s)}Y)\\
=&-(\overline{\nabla}_{\overline{\nabla} r}\overline{\hess} r)(Y)\\
=&-(\overline{\nabla}_Y\overline{\hess} r)(\overline{\nabla} r)+\overline{\R}(\overline{\nabla} r,Y)\overline{\nabla} r\\
=&\overline{\hess} r(\overline{\nabla}_Y\overline{\nabla} r)+\overline{\R}(\overline{\nabla} r,Y)\overline{\nabla} r.
\end{align*}
On the other hand, since $Y$ is parallel
\begin{displaymath}
\frac{d}{dt}\pair{\overline{\hess} r(Y),Y}\Big|_s=\pair{\overline{\nabla}_{\gamma'(s)}\overline{\hess} r(Y),Y}.
\end{displaymath}
Hence
\begin{displaymath}
\frac{d}{dt}\overline{\Hess} r(\gamma)(Y,Y)-\pair{\overline{\hess} r(\gamma)(Y),\overline{\hess} r(\gamma)(Y)}=-{K_M}_{\gamma}(Y\wedge\gamma')
\end{displaymath}
Notice that
\begin{displaymath}
\overline{\Hess} r(X,X)\geq\lambda_{\min}
\end{displaymath}
for every spacelike unit vector field $X\bot\overline{\nabla} r$. Let us choose $Y$ so that at $s$
\begin{displaymath}
\overline{\Hess} r(\gamma)(Y,Y)=\lambda_{\min}(\gamma(s)).
\end{displaymath}
Then, the function $\Hess r(\gamma)(Y,Y)-\lambda_{\min}\circ\gamma$ attains its minimum at $s$. Hence
\begin{displaymath}
\frac{d}{dt}\overline{\Hess} r(\gamma)(Y,Y)\Big|_s=\frac{d}{dt}(\lambda_{\min}\circ\gamma)\Big|_s
\end{displaymath}
and we have proved that $\lambda_{\min}$ satisfies the first equation in \eqref{lambdamin}, since $K(Y\wedge\gamma')\leq G$.
The asymptotic behaviour follows from the expression
\begin{equation}\label{hessnear0}
\overline{\Hess} r=\frac1r(\pair{,}+dr\otimes dr)+o(1)
\end{equation}
that can be proved using normal coordinates around $p$. Now, if we set $\phi=\frac{h'}{h}$, we find that $\phi$ satisfies
\begin{displaymath}
\left\{\begin{array}{l}
\phi'+\phi^2=G\qquad\text{on }(0,r_0)\\
\phi=\frac1t+o(t)\qquad\text{as }t\ra0^+
\end{array}
\right.
\end{displaymath}
Then, using Corollary \ref{riccati} with $g_1=\lambda_{\min}$, $g_2=\phi$ and $\alpha=1$ we conclude that
\[
\lambda_{\min}(q)\geq -\frac{h'}{h}(r(q)) 
\]
and this concludes the proof of $(a)$.\\
Finally, for what concerns claim $(b)$, we observe that reasoning as in the proof
of claim $(a)$ and choosing $Y$ so that at $s$
\begin{displaymath}
\overline{\Hess} r(\gamma)(Y,Y)=\lambda_{\max}(\gamma(s))
\end{displaymath}
we can prove that, if \eqref{curvlowbound} holds, 
$\lambda_{\max}$ satisfies
\begin{displaymath}
\left\{\begin{array}{l}
\frac{d}{dt}(\lambda_{\max}\circ \gamma)-(\lambda_{\max}\circ\gamma)^2\leq- G\qquad\text{for a.e. }t>0\\
\lambda_{\max}\circ\gamma=\frac1t+o(t)\qquad\text{as }t\ra0^+
\end{array}
\right.
\end{displaymath}
In this case, setting again $\phi=\frac{h'}{h}$, we find that $\phi$ satisfies
\begin{displaymath}
\left\{\begin{array}{l}
\phi'+\phi^2=G\qquad\text{on }(0,r_0)\\
\phi=\frac1t+o(t)\qquad\text{as }t\ra0^+
\end{array}
\right.
\end{displaymath}
Then, we can conclude again using Corollary \ref{riccati} with $g_1=\lambda_{\max}$, $g_2=\phi$ and $\alpha=1$.
\end{proof}
\begin{theorem}[Laplacian Comparison Theorem]
Let $M^{n+1}$ be an $n+1$- dimensional spacetime. Assume that there exists a point $p \in M$ such that $\mathcal{I}^+(p)\neq\emptyset$ and let $q \in \mathcal{I}^+(p)$. Let $r(\cdot)=d_p(\cdot)$ be the Lorentzian distance function from $p$. Given a smooth even function $G$ on $\erre$, let $h$ be a solution of the Cauchy problem
\begin{displaymath}
\left\{\begin{array}{l}
h''-Gh\geq0\\
h(0)=0,\ h'(0)=1
\end{array}\right.
\end{displaymath}
and let $I=[0,r_0)\subset[0,+\infty)$ be the maximal interval where $h$ is positive. If
\begin{equation}
\ricc_{M}(\overline{\nabla} r,\overline{\nabla} r)\geq -nG(r),
\end{equation}
then
\begin{displaymath}
\overline{\Delta} r\geq -n\frac{h'}{h}(r)
\end{displaymath}
holds pointwise on $\mathcal{I}^+(p)\cap B^+(p,r_0)$.
\end{theorem}
\begin{proof}
Let $v\in\exp_p^{-1}(q)\in\rm{int}(\tilde{\mathcal{I}}^+(p))$ and let $\gamma(t)=\exp_p(tv)$, $0\leq t\leq s_p(v)$, be the radial future directed unit timelike geodesic with $\gamma(0)=p$, $\gamma(s)=q$, $s=r(q)$. Recall that $\gamma'(s)=-\overline{\nabla} r(q)$ and $\overline{\nabla}_{\overline{\nabla} r}\overline{\nabla} r(q)=0$. Define
\begin{displaymath}
\varphi(t)=\overline{\Delta} r\circ\gamma(t),\qquad t\in(0,s].
\end{displaymath}
Then tracing Equation \eqref{hessnear0}
\begin{displaymath}
\varphi(t)=\frac nt+o(t)\qquad\text{as }t\ra0^+.
\end{displaymath}
Recall that given $f\in C^{\infty}(M)$ the following Bochner formula holds
$$
\frac{1}{2}\overline{\Delta} \pair{\overline{\nabla} f,\overline{\nabla} f}=\norm{\overline{\hess} f}^2+\ricc_M(\overline{\nabla} f,\overline{\nabla} f)+\pair{\overline{\nabla} \overline{\Delta} f,\overline{\nabla} f}. 
$$
See \cite{grk} for more details.
Since $\norm{\overline{\nabla} r}^2=-1$, it follows that
\begin{displaymath}
0=\norm{\overline{\hess} r}^2+\ricc_M(\overline{\nabla} r,\overline{\nabla} r)+\pair{\overline{\nabla} \overline{\Delta} r,\overline{\nabla} r}.
\end{displaymath}
Since $\norm{\overline{\hess} r}^2\geq\frac{(\overline{\Delta} r)^2}{n}$ and $\ricc_M(\overline{\nabla} r,\overline{\nabla} r)\geq -nG(r)$, we have
\begin{displaymath}
\frac{1}{n}(\overline{\Delta} r)^2+\pair{\overline{\nabla} \overline{\Delta}  r,\overline{\nabla} r}\leq nG(r).
\end{displaymath}
Computing along $\gamma$
\begin{displaymath}
\varphi'(t)=\frac{d}{dt}(\overline{\Delta} r(\gamma(t)))\Big|_s=\pair{\overline{\nabla} \overline{\Delta} r(\gamma(t)),\gamma'(t)}\Big|_s=-\pair{\overline{\nabla} \overline{\Delta} r,\overline{\nabla} r}.
\end{displaymath}
Hence the function $\varphi$ satisfies
\begin{displaymath}
\left\{\begin{array}{l}
\varphi'(t)-\frac{\varphi^2(t)}{n}\geq -n G\\
\varphi(t)=\frac nt +o(t)\qquad\text{as }t\ra0^+
\end{array}
\right.
\end{displaymath}
Set $\phi=n\frac{h'}{h}$. Then $\phi$ satisfies
\begin{displaymath}
\left\{\begin{array}{l}
\phi'(t)+\frac{\phi^2(t)}{n}\geq n G\qquad\text{on }(0,r_0)\\
\phi(t)=\frac nt +o(t)\qquad\text{as }t\ra0^+
\end{array}
\right.
\end{displaymath}
Then we conclude again using Corollary \ref{riccati}.
\end{proof}
\section{Applications to spacelike hypersurfaces}

Let $\psi:\Sigma^n\ra M^{n+1}$ be a spacelike hypersurface isometrically immersed into the spacetime $M$. Since $M$ is time-orientable, there exists a unique future-directed timelike unit normal field $\nu$ globally defined on $\Sigma$. We will refer to that normal field $\nu$ as the future-pointing Gauss map of the hypersurface.
We let $A: T\Sigma \rightarrow T\Sigma $ denote the second fundamental form of the immersion. Its eigenvalues $k_1,...,k_n$ are the principal curvatures of the hypersurface. Their elementary symmetric functions 
\begin{align*}
S_k=&\sum_{i_1<...<i_k}k_{i_1}\cdots k_{i_k},\qquad k=1,...,n,\\ S_0=&1,
\end{align*} 
define the $k$-mean curvatures of the immersion via the formula
$$
{n \choose k}H_k= (-1)^kS_k.
$$ 
Thus $H_1=-1/n\mathrm{Tr}(A)=H$ is the mean curvature of $\Sigma$ and  $n(n-1)H_2=\overline{S}-S+2\overline{\ricc}(\nu,\nu)$, where $S$ and $\overline{S}$ are, respectively, the scalar curvature of $\Sigma$ and $M^{n+1}$ and $\overline{\ricc}$ is the Ricci tensor of $M^{n+1}$. Even more, when $k$ is even, it follows from the Gauss equation that $H_k$ is a geometric quantity which is related to the intrinsic curvature of $\Sigma^n$.

The classical Newton transformations associated to the immersion are defined inductively by
$$
P_0=I,\qquad P_k={n \choose k}H_kI+AP_{k-1},
$$
for every $k=1,...,n$.
\begin{prop} The following formulas hold:
\begin{enumerate}
\item $\mathrm{Tr}(P_k)=c_k H_k$,\\
\item $\mathrm{Tr}(AP_k)=-c_k H_{k+1}$,\\
\item $\mathrm{Tr}(A^2P_k)={n\choose {k+1}}(nH_1H_{k+1}-(n-k-1)H_{k+2})$,
\end{enumerate}
where $c_k=(n-k){n \choose k}=(k+1){n\choose {k+1}}$.
\end{prop}
We refer the reader to \cite{aliasbrasilcolares} for the proof of the last proposition and for further details on the Newton transformations (see also \cite{reilly} and \cite{rosenberg} for others details on the Newton transformations in the Riemannian setting). 
Let $\nabla$ be the Levi-Civita connection of $\Sigma$. We define the second order linear differential operator $L_k:C^{\infty}(\Sigma) \ra C^{\infty}(\Sigma)$ associated to $P_k$ by
$$
L_k f=\mathrm{Tr}(P_k \circ \hess f).
$$
It follows by the definition that the operator $L_k$ is elliptic if and only if $P_k$ is positive definite. Let us state two useful 
lemmas in which geometric conditions are given in order to guarantee the ellipticity of $L_k$ when $k \geq 1$ (Recall that 
$L_0=\Delta$ is always elliptic).
\begin{lemma}\label{lemmaelliptl1lor}
Let $\Sigma$ be a spacelike hypersurface immersed into a spacetime. If $H_2 >0$ on $\Sigma$, then $L_1$ is an elliptic operator 
(for an appropriate choice of the Gauss map $\nu$).
\end{lemma}
For a proof of Lemma \ref{lemmaelliptl1lor} see Lemma 3.10 in \cite{elbert}. The next Lemma is a consequence of Proposition 3.2 in 
\cite{barbosacolares}.
\begin{lemma}\label{ellipticitylr}
Let $\Sigma^n$ be a spacelike hypersurface immersed into a $n+1$-dimensional spacetime. If there exists an elliptic point of 
$\Sigma$, with respect to an appropriate choice of the Gauss map $\nu$, and $H_k >0$ on $\Sigma$, $3 \leq k \leq n$, then for all 
$1 \leq j \leq k-1$ the operator $L_j$ is elliptic.
\end{lemma}
We recall here that given a spacelike hypersurface $\Sigma$, a point $p\in\Sigma$ is said to be elliptic if the second fundamental 
form of the immersion is negative definite at $p$.

Now consider $\psi:\Sigma^n\ra M^{n+1}$ and assume that there exists a point $p \in M$ such that $\mathcal{I}^+(p)\neq \emptyset$ and that $\psi(\Sigma)\subset\mathcal{I}^+(p)$. Let $r(\cdot)=d_p(\cdot)$ be the Lorentzian distance function from $p$ and let $u=r\circ \psi:\Sigma\ra(0,+\infty)$ be the function $r$ along the hypersurface, which is a smooth function on $\Sigma$. Let us calculate the Hessian of $u$ on $\Sigma$. Notice that
\begin{displaymath}
\overline{\nabla}r=\nabla u-\pair{\overline{\nabla}r,\nu}\nu. 
\end{displaymath}
Hence, since $\norm{\overline{\nabla}r}^2=-1$ and $\pair{\overline{\nabla}r,\nu}>0$, we have
\begin{displaymath}
\pair{\overline{\nabla}r,\nu}=\sqrt{1+\norm{\nabla u}^2}\geq 1.
\end{displaymath}
Hence
\begin{displaymath}
\overline{\nabla}r=\nabla u-\nu\sqrt{1+\norm{\nabla u}^2}
\end{displaymath}
Moreover
\begin{displaymath}
\overline{\nabla}_X\overline{\nabla}r=\nabla_X\nabla u+\sqrt{1+\norm{\nabla u}^2}AX+\pair{AX,\nabla u}\nu-X(\sqrt{1+\norm{\nabla u}^2})\nu
\end{displaymath}
for every spacelike $X \in T\Sigma$. Thus
\begin{displaymath}
\Hess u(X,P_kX)=\overline{\Hess} r(X,P_k X)-\sqrt{1+\norm{\nabla u}^2}\pair{P_kA X,X}
\end{displaymath} 
On the other hand, we have the following decompositions
\begin{align*}
X=&X^*-\pair{X,\nabla u}\overline{\nabla} r\\
P_kX=&(P_kX)^*-\pair{X,P_k\nabla u}\overline{\nabla} r,
\end{align*}
where $X^*$, $(P_kX)^*$ are respectively the components of $X$, $P_kX$ orthogonal to $\overline{\nabla}r$.
Then
\begin{displaymath}
\pair{X^*,(P_kX)^*}=\pair{X,P_kX}+\pair{X,P_k\nabla u}\pair{X,\nabla u}
\end{displaymath}
and, taking into account that
\begin{displaymath}
\overline{\nabla}_{\overline{\nabla}r}\overline{\nabla}r=0
\end{displaymath}
we find
\begin{displaymath}
\overline{\Hess}r(X,P_kX)=\overline{\Hess}r(X^*,(P_kX)^*).
\end{displaymath}
Hence, if we assume that $K_M(\Pi)\leq G(r)$ for all timelike planes $\Pi$, then
\begin{align*}
\overline{\Hess}r(X,P_kX)=&\overline{\Hess}r(X^*,(P_kX)^*)\geq -\frac{h'}{h}(u)\pair{X^*,(P_k X)^*}\\
=&-\frac{h'}{h}(u)(\pair{X,P_kX}+\pair{X,\nabla u}\pair{X,P_k\nabla u}),
\end{align*}
where $h$ is a solution of the problem
\begin{displaymath}
\left\{\begin{array}{l}
h''-Gh=0\\
h(0)=0,\ h'(0)=1
\end{array}
\right.
\end{displaymath}
Therefore
\begin{align*}
\Hess u(X,P_kX)\geq &-\frac{h'}{h}(u)(\pair{X,P_kX}+\pair{X,\nabla u}\pair{X,P_k\nabla u})\\
&-\sqrt{1+\norm{\nabla u}^2}\pair{P_kA X,X}.
\end{align*} 
Tracing
\begin{displaymath}
L_ku\geq-\frac{h'}{h}(u)(c_kH_k+\pair{\nabla u,P_k\nabla u})+\sqrt{1+\norm{\nabla u}^2}c_kH_{k+1}.
\end{displaymath} 
Summarizing, we have proved the following
\begin{prop}\label{proplkugreater}
Let $M^{n+1}$ be an $n+1$-dimensional spacetime. Assume that there exists a point $p \in M$ such that $\mathcal{I}^+(p)\neq\emptyset$ 
and let $r(\cdot)=d_p(\cdot)$ be the Lorentzian distance function from $p$. Given a smooth even function $G$ on $\erre$, let $h$ be a 
solution of the Cauchy problem
\begin{displaymath}
\left\{\begin{array}{l}
h''-Gh= 0\\
h(0)=0,\ h'(0)=1
\end{array}\right.
\end{displaymath}
and let $I=[0,r_0)\subset[0,+\infty)$ be the maximal interval where $h$ is positive. Let $\psi:\Sigma^n\ra M^{n+1}$ be a spacelike 
hypersurface such that $\psi(\Sigma^n)\subset\mathcal{I}^+(p)\cap B^+(p,r_0)$. If
\begin{equation}
K_{M}(\Pi)\leq G(r)
\end{equation}
for all timelike planes $\Pi$, then
\begin{equation}
L_ku\geq-\frac{h'}{h}(u)(c_kH_k+\pair{\nabla u,P_k\nabla u})+\sqrt{1+\norm{\nabla u}^2}c_kH_{k+1}.
\end{equation}
\end{prop}
On the other hand, if we assume that $K_M(\Pi)\geq G(r)$ for all timelike planes in $M$, the same computations yield the following
\begin{prop}\label{proplkuless}
Let $M^{n+1}$ be an $n+1$-dimensional spacetime. Assume that there exists a point $p \in M$ such that $\mathcal{I}^+(p)\neq\emptyset$ 
and let $r(\cdot)=d_p(\cdot)$ be the Lorentzian distance function from $p$. Given a smooth even function $G$ on $\erre$, let $h$ be a 
solution of the Cauchy problem
\begin{displaymath}
\left\{\begin{array}{l}
h''-Gh= 0\\
h(0)=0,\ h'(0)=1
\end{array}\right.
\end{displaymath}
and let $I=[0,r_0)\subset[0,+\infty)$ be the maximal interval where $h$ is positive. Let $\psi:\Sigma^n\ra M^{n+1}$ be a spacelike 
hypersurface such that $\psi(\Sigma^n)\subset\mathcal{I}^+(p)\cap B^+(p,r_0)$. If
\begin{equation}
K_{M}(\Pi)\geq G(r)
\end{equation}
for all timelike planes $\Pi$, then
\begin{equation}
L_ku\leq-\frac{h'}{h}(u)(c_kH_k+\pair{\nabla u,P_k\nabla u})+\sqrt{1+\norm{\nabla u}^2}c_kH_{k+1}.
\end{equation}
\end{prop}
In the following, under suitable bounds on the sectional curvature of the ambient spacetime, we will find some lower and upper bounds 
for the mean cuevature and the higher order mean curvatures associated to the immersion. In order to do it we will use the 
\textit{Omori-Yau maximum principle} for the Laplacian and for more general elliptic operators (for more details and others 
applications of this technique see \cite{aliasimperarigoli}, \cite{aliasimperarigolilor}). Namely, if $L=\text{Tr}(P\circ\hess)$, 
where $P$ is a symmetric operator with trace bounded above, using the terminology introduced by S. Pigola, M. Rigoli and A. G. Setti 
in \cite{pirise},we say that the \textit{ 
Omori-Yau maximum principle} holds on $\Sigma$ for $L$ if
for any smooth function $u\in C^{\infty}(\Sigma)$ with $u^*=\sup_{\Sigma}u<+\infty$ there exists a sequence of points $\{p_i\}_{i\in \enne}\subset \Sigma$ such that
\begin{equation}\label{genomoriyau}
(i)\ u(p_i)>u^*-\frac{1}{i},\ (ii) \ \norm{\nabla u(p_i)}<\frac{1}{i}, \ (iii) \ Lu(p_i)< \frac{1}{i}.
\end{equation}
Equivalently if $u_*=\inf_{\Sigma} u>-\infty$, we can find a sequence $\{q_i\}_{i \in\enne}\subset\Sigma$ such that
\begin{equation}\label{genomoriyau}
(i)\ u(q_i)>u_*-\frac{1}{i},\ (ii) \ \norm{\nabla u(q_i)}<\frac{1}{i}, \ (iii) \ Lu(q_i)> -\frac{1}{i}.
\end{equation}
Clearly the Laplacian belong to this class of operators. In this case, S. Pigola, M. Rigoli and A. G. Setti showed in \cite{pirise} 
that a condition of the form
\begin{equation}\label{condricci}
\ricc(\nabla\rho,\nabla\rho)\geq-C^2G(\rho),
\end{equation}
where $\rho$ is the distance function on $\Sigma$ to a fixed point and $G:[0,+\infty)\ra\erre$ is a smooth function satisfying
\begin{equation}\label{condG}
\begin{array}{ll}
(i)\  G(0)>0, & (ii)\  G'(t)\geq 0 \qquad \text{on } [0,+\infty),\\
(iii)\  G(t)^{-\frac{1}{2}}\not \in L^1(+\infty),& (iv)\  \limsup_{t \ra \infty} \frac{t G(\sqrt{t})}{G(t)}<+\infty.  
\end{array}
\end{equation}
is sufficient to guarantee the validity of the Omori-Yau maximum principle for the Laplacian on $\Sigma$.
Analogously, in \cite{aliasimperarigoli}, L. J. Alias, M. Rigoli and the author showed that the condition 
\begin{equation}\label{condseccurv}
K(\nabla \rho,X)\geq -G(\rho),
\end{equation}
where $X$ is any vector field tangent to $\Sigma$ and $G$ satisfies \eqref{condG}, is sufficient to guarantee the validity of the 
Omori-Yau maximum principle on $\Sigma$ for operators $L$ with the properties described above. \\
Applying the Omori-Yau maximum principle we find the following estimates for the mean curvature. The proof of the following theorems 
is essentially the same as that of Theorems 4.1 and 4.2 in \cite{aliashurtadopalmer}.\\
\begin{theorem}
Let $M^{n+1}$ be an $n+1$-dimensional spacetime. Assume that there exists a point $p \in M$ such that $\mathcal{I}^+(p)\neq\emptyset$ and let $r(\cdot)=d_p(\cdot)$ be the Lorentzian distance function from $p$. Given a smooth even function $G$ on $\erre$, let $h$ be a solution of the Cauchy problem
\begin{displaymath}
\left\{\begin{array}{l}
h''-Gh= 0\\
h(0)=0,\ h'(0)=1
\end{array}\right.
\end{displaymath}
and let $I=[0,r_0)\subset[0,+\infty)$ be the maximal interval where $h$ is positive. Let $\psi:\Sigma^n\ra M^{n+1}$ be a spacelike 
hypersurface such that $\psi(\Sigma^n)\subset\mathcal{I}^+(p)\cap B^+(p,\delta)$ with $\delta\leq r_0$. If
\begin{equation}
\ricc_{M}(\overline{\nabla} r, \overline{\nabla} r)\geq -nG(r),
\end{equation}
and the Omori-Yau maximum principle holds on $\Sigma$, then 
\begin{displaymath}
\inf_{\Sigma}H_1\leq\frac{h'}{h}\big(\sup_{\Sigma}u\big),
\end{displaymath}
where $u$ denotes the Lorentzian distance $d_p$ along the hypersurface. 
\end{theorem}
On the other hand, if we assume that the sectional curvature of timelike planes is bounded from below we obtain
\begin{theorem}\label{lowboundH1}
Let $M^{n+1}$ be an $n+1$- dimensional spacetime. Assume that there exists a point $p \in M$ such that 
$\mathcal{I}^+(p)\neq\emptyset$ and let $r(\cdot)=d_p(\cdot)$ be the Lorentzian distance function from $p$. Given a smooth even 
function $G$ on $\erre$, let $h$ be a solution of the Cauchy problem
\begin{displaymath}
\left\{\begin{array}{l}
h''-Gh= 0\\
h(0)=0,\ h'(0)=1
\end{array}\right.
\end{displaymath}
and let $I=[0,r_0)\subset[0,+\infty)$ be the maximal interval where $h$ is positive and let $\psi:\Sigma^n\ra M^{n+1}$ be a 
spacelike hypersurface such that $\psi(\Sigma^n)\subset\mathcal{I}^+(p)\cap B^+(p,r_0)$. If
\begin{equation}
K_{M}(\Pi)\geq G(r)
\end{equation}
for all timelike planes $\Pi$ and if the Omori-Yau maximum principle holds on $\Sigma$, then
\begin{displaymath}
\sup_{\Sigma}H_1\geq\frac{h'}{h}\big(\inf_{\Sigma}u\big),
\end{displaymath}
where $u$ denotes the Lorentzian distance $d_p$ along the hypersurface $\Sigma$. 
\end{theorem}

The previous estimates can be extended to the higher order mean curvatures in the following way. To find the estimates we will use 
the Omori-Yau maximum principle for elliptic operators of the form $L=\text{Tr}(P\circ\hess)$, where $P$ is a symmetric operator 
with trace bounded above. For simplicity, we will refer to that as the generalized Omori-Yau maximum principle.
\begin{theorem}
Let $M^{n+1}$ be an $n+1$-dimensional spacetime. Assume that there exists a point $p \in M$ such that $\mathcal{I}^+(p)\neq\emptyset$ 
and let $r(\cdot)=d_p(\cdot)$ be the Lorentzian distance function from $p$. Given a smooth even function $G$ on $\erre$, let $h$ be a 
solution of the Cauchy problem
\begin{displaymath}
\left\{\begin{array}{l}
h''-Gh= 0\\
h(0)=0,\ h'(0)=1
\end{array}\right.
\end{displaymath}
and let $I=[0,r_0)\subset[0,+\infty)$ be the maximal interval where $h$ is positive. Let $\psi:\Sigma^n\ra M^{n+1}$ be a spacelike 
hypersurface such that $\psi(\Sigma^n)\subset\mathcal{I}^+(p)\cap B^+(p,\delta)$, with $\delta\leq r_0$. Assume that $H_2>0$ and that 
$\sup_\Sigma H_1<+\infty$. If
\begin{equation}
K_{M}(\Pi)\leq G(r)
\end{equation}
for all timelike planes $\Pi$ and if Omori-Yau maximum principle holds on $\Sigma$, then 
\begin{displaymath}
\inf_{\Sigma}H_2^{\frac 12}\leq\Big|\frac{h'}{h}\big(\sup_{\Sigma}u\big)\Big|,
\end{displaymath}
where $u$ denotes the Lorentzian distance $d_p$ along the hypersurface. 
\end{theorem}
\begin{proof}
Consider the operator
\begin{align*}
\mathcal{L}=&L_1 +(n-1)\frac{1}{\sqrt{1+\norm{\nabla u}^2}}\Big(\Big|\frac{h'}{h}(u)\Big|\Big)\Delta\\
=& \text{Tr}(\mathcal{P}\circ\hess),
\end{align*}
where
\begin{displaymath}
\mathcal{P}=P_1+(n-1)\frac{1}{\sqrt{1+\norm{\nabla u}^2}}\Big(\Big|\frac{h'}{h}(u)\Big|\Big)I.
\end{displaymath}
Notice that, since $H_2>0$, the operator $L_1$ is elliptic and so is $\mathcal{L}$. Since $0<u<\sup_{\Sigma}u<\delta$, $h'/h(u)$ is bounded. Furthermore, $\sup_{\Sigma}H_1<+\infty$ and $1/\sqrt{1+\norm{\nabla u}^2}\leq1$, hence we can apply the Omori-Yau maximum principle for the operator $\mathcal{L}$. We can then find a sequence $\{p_i\}_{i\in\enne}\subset\Sigma$ such that
\begin{displaymath}
(i)\ u(p_i)>u^*-\frac{1}{i},\ (ii) \ \norm{\nabla u(p_i)}<\frac{1}{i}, \ (iii) \ \mathcal{L}u(p_i)< \frac{1}{i}.
\end{displaymath}
A simple computation using Proposition \ref{proplkugreater} shows that
\begin{align*}
\mathcal{L}u\geq& -(n-1)\frac{1}{\sqrt{1+\norm{\nabla u}^2}}\Big(\frac{h'}{h}(u)\Big)^2(n+\norm{\nabla u}^2)-\\
&-\Big(\Big|\frac{h'}{h}(u)\Big|\Big)\pair{P_1 \nabla u,\nabla u}+n(n-1)\sqrt{1+\norm{\nabla u}^2}H_2\\
\geq& -(n-1)\frac{1}{\sqrt{1+\norm{\nabla u}^2}}\Big(\frac{h'}{h}(u)\Big)^2(n+\norm{\nabla u}^2)-\\
&-\Big(\Big|\frac{h'}{h}(u)\Big|\Big)\pair{P_1 \nabla u,\nabla u}+n(n-1)\sqrt{1+\norm{\nabla u}^2}\inf_\Sigma H_2.
\end{align*}
Hence
\begin{align*}
\frac{1}{i}>\mathcal{L}u(p_i)\geq & -(n-1)\frac{1}{\sqrt{1+\norm{\nabla u(p_i)}^2}}\Big(\frac{h'}{h}(u(p_i))\Big)^2(n+\norm{\nabla u(p_i)}^2)-\\
&-\Big(\Big|\frac{h'}{h}(u(p_i))\Big|\Big)\pair{P_1 \nabla u(p_i),\nabla u(p_i)}\\
&+n(n-1)\sqrt{1+\norm{\nabla u(p_i)}^2}\inf_\Sigma H_2.
\end{align*}
Taking the limit for $i\ra+\infty$ we find
\begin{displaymath}
0\geq -n(n-1)\Big(\frac{h'}{h}(\sup_{\Sigma}u)\Big)^2+n(n-1)\inf_{\Sigma}H_2.
\end{displaymath}
and the conclusion follows.
\end{proof}
\begin{theorem}\label{boundinfHk}
Let $M^{n+1}$ be an $n+1$- dimensional spacetime, $n \geq 3$. Assume that there exists a point $p \in M$ such that $\mathcal{I}^+(p)\neq\emptyset$ and let $r(\cdot)=d_p(\cdot)$ be the Lorentzian distance function from $p$. Given a smooth even function $G$ on $\erre$, let $h$ be a solution of the Cauchy problem
\begin{displaymath}
\left\{\begin{array}{l}
h''-Gh= 0\\
h(0)=0,\ h'(0)=1
\end{array}\right.
\end{displaymath}
and let $I=[0,r_0)\subset[0,+\infty)$ be the maximal interval where $h$ is positive. Let $\psi:\Sigma^n\ra M^{n+1}$ be a spacelike 
hypersurface such that $\psi(\Sigma^n)\subset\mathcal{I}^+(p)\cap B^+(p,\delta)$, with $\delta\leq r_0$. 
Assume that there exists an elliptic point $p_0\in\Sigma$, that $H_k>0$, $3\leq k\leq n$, and that $\sup_{\Sigma}H_1<+\infty$. If
\begin{equation}
K_{M}(\Pi)\leq G(r)
\end{equation}
for all timelike planes $\Pi$ and if the generalized Omori-Yau maximum principle holds on $\Sigma$, then 
\begin{displaymath}
\inf_{\Sigma}H_k^{\frac 1k}\leq\Big|\frac{h'}{h}\big(\sup_{\Sigma}u\big)\Big|,
\end{displaymath}
where $u$ denotes the Lorentzian distance $d_p$ along the hypersurface. 
\end{theorem}
\begin{proof}
Consider the operator
\begin{displaymath}
\mathcal{L}=\sum_{j=0}^{k-1}(1+\norm{\nabla u}^2)^{-\frac {k-1-j}{2}}\Big(\Big|\frac{h'}{h}(u)\Big|\Big)^{k-1-j}\frac{c_{k-1}}{c_j}L_j
\end{displaymath}
Notice that, since there exists an elliptic point $p_0\in\Sigma$ and $H_k>0$, $3\leq k\leq n$, the operators $L_j$ are elliptic for 
all $1\leq j\leq k-1$. Since $0<u<\sup_{\Sigma}u<\delta$, $1/\sqrt{1+\norm{\nabla u}^2}\leq1$ and $\sup_{\Sigma}H_1<+\infty$, we can 
apply the Omori-Yau maximum principle for the operator $\mathcal{L}$. Hence, we can find a sequence $\{p_i\}_{i\in\enne}\subset\Sigma$ 
such that
\begin{displaymath}
(i)\ u(p_i)>u^*-\frac{1}{i},\ (ii) \ \norm{\nabla u(p_i)}<\frac{1}{i}, \ (iii) \ \mathcal{L}u(p_i)< \frac{1}{i}.
\end{displaymath}
A straightforward computation using Proposition \ref{proplkugreater} shows that
\begin{align*}
\mathcal{L}u\geq& -\sum_{j=0}^{k-1}(1+\norm{\nabla u}^2)^{-\frac {k-1-j}{2}}\Big(\Big|\frac{h'}{h}(u)\Big|\Big)^{k-j}\frac{c_{k-1}}{c_j}\pair{P_j\nabla u,\nabla u}\\
&-c_{k-1}\frac{1}{(1+\norm{\nabla u}^2)^{(k-1)/2}}\Big(\Big|\frac{h'}{h}(u)\Big|\Big)^k+\sqrt{1+\norm{\nabla u}^2}c_{k-1}H_k.
\end{align*}
Hence
\begin{align*}
\frac{1}{i}>\mathcal{L}u(p_i)\geq &-c_{k-1}\frac{1}{(1+\norm{\nabla u(p_i)}^2)^(k-1)/2}\Big(\Big|\frac{h'}{h}(u(p_i))\Big|\Big)^k\\ &-\sum_{j=1}^{k-1}(1+\norm{\nabla u(p_i)}^2)^{-\frac {k-1-j}{2}}\Big(\Big|\frac{h'}{h}(u(p_i))\Big|\Big)^{k-j}\frac{c_{k-1}}{c_j}\pair{P_j\nabla u,\nabla u}(p_i)\\
&+\sqrt{1+\norm{\nabla u(p_i)}^2}c_{k-1}\inf_\Sigma H_k. 
\end{align*}
Taking the limit for $i\ra+\infty$ we find
\begin{displaymath}
0\geq -c_{k-1}\Big(\Big|\frac{h'}{h}(\sup_{\Sigma}u)\Big|\Big)^k+c_{k-1}\inf_{\Sigma}H_k.
\end{displaymath}
\end{proof}
On the other hand, if we assume that the sectional curvature of timelike planes is bounded from below we find the following estimates
\begin{theorem}\label{boundsupH2}
Let $M^{n+1}$ be an $n+1$-dimensional spacetime. Assume that there exists a point $p \in M$ such that $\mathcal{I}^+(p)\neq\emptyset$ 
and let $r(\cdot)=d_p(\cdot)$ be the Lorentzian distance function from $p$. Given a smooth even function $G$ on $\erre$, let $h$ be a 
solution of the Cauchy problem
\begin{displaymath}
\left\{\begin{array}{l}
h''-Gh= 0\\
h(0)=0,\ h'(0)=1
\end{array}\right.
\end{displaymath}
and let $I=[0,r_0)\subset[0,+\infty)$ be the maximal interval where $h$ is positive. Let $\psi:\Sigma^n\ra M^{n+1}$ be a spacelike 
hypersurface such that $\psi(\Sigma^n)\subset\mathcal{I}^+(p)\cap B^+(p,r_0)$. Assume that $H_2>0$ and that $\sup_{\Sigma}H_1<+\infty$. If
\begin{equation}
K_{M}(\Pi)\geq G(r)
\end{equation}
for all timelike planes $\Pi$ and if the generalized Omori-Yau maximum principle holds on $\Sigma$, then 
\begin{displaymath}
\sup_{\Sigma}H_2^{\frac 12}\geq\frac{h'}{h}\big(\inf_{\Sigma}u\big),
\end{displaymath}
where $u$ denotes the Lorentzian distance $d_p$ along the hypersurface. 
\end{theorem}
\begin{theorem}\label{boundsupHk}
Let $M^{n+1}$ be an $n+1$- dimensional spacetime, $n \geq 3$. Assume that there exists a point $p \in M$ such that $\mathcal{I}^+(p)\neq\emptyset$ and let $r(\cdot)=d_p(\cdot)$ be the Lorentzian distance function from $p$. Given a smooth even function $G$ on $\erre$, let $h$ be a solution of the Cauchy problem
\begin{displaymath}
\left\{\begin{array}{l}
h''-Gh= 0\\
h(0)=0,\ h'(0)=1
\end{array}\right.
\end{displaymath}
and let $I=[0,r_0)\subset[0,+\infty)$ be the maximal interval where $h$ is positive. Let $\psi:\Sigma^n\ra M^{n+1}$ be a spacelike 
hypersurface such that $\psi(\Sigma^n)\subset\mathcal{I}^+(p)\cap B^+(p,r_0)$. Assume that there exists an elliptic point 
$p_0\in\Sigma$, that $H_k>0$, $3\leq k\leq n$, and that $\sup_{\Sigma}H_1<+\infty$. If
\begin{equation}
K_{M}(\Pi)\geq G(r)
\end{equation}
for all timelike planes $\Pi$ and if the generalized Omori-Yau maximum principle holds on $\Sigma$, then
\begin{displaymath}
\sup_{\Sigma}H_k^{\frac 1k}\geq\frac{h'}{h}\big(\inf_{\Sigma}u\big),
\end{displaymath}
where $u$ denotes the Lorentzian distance $d_p$ along the hypersurface. 
\end{theorem}
We will only prove Theorem \ref{boundsupHk}. The proof of Theorem \ref{boundsupH2} proceed exactly in the same way.   
\begin{proof}[Proof of Theorem \ref{boundsupHk}]
If $h'/h(\inf_{\Sigma} u)\leq0$, the result is trivial since
\begin{displaymath}
\frac{h'}{h}\Big(\inf_{\Sigma}u\Big)\leq 0< \sup_{\Sigma}H_k^{\frac 1k}.
\end{displaymath}
Conversely, assume $h'/h(\inf_{\Sigma}u)>0$. Since $u\geq u_*:=\inf_\Sigma u\geq 0$, we want to apply the Omori-Yau maximum principle 
for a suitable elliptic operator with trace bounded above. Notice that it must be $\inf_{\Sigma} u>0$. Indeed, if $\inf_{\Sigma} u=0$, 
since $\lim_{s\ra0}h'/h(s)=+\infty$, it follows by the estimate in Theorem \ref{lowboundH1} that $\sup_{\Sigma}H_1=+\infty$, which 
contradicts our assumptions.  The operator that we consider is the following
\begin{align*}
\mathcal{L}=&\sum_{j=0}^{k-1}(1+\norm{\nabla u}^2)^{-\frac{k-j-1}{2}}\Big(\frac{h'}{h}(\inf_{\Sigma}u)\Big)^{k-j-1}\frac{c_{k-1}}{c_j}L_j\\
=&\text{Tr}(\mathcal{P}\circ\hess),
\end{align*}
where 
$$
\mathcal{P}=\sum_{j=0}^{k-1}(1+\norm{\nabla u}^2)^{-\frac{k-j-1}{2}}\Big(\frac{h'}{h}(\inf_{\Sigma}u)\Big)^{k-j-1}\frac{c_{k-1}}{c_j}P_j.
$$
Notice that, since there exists an elliptic point $p_0\in\Sigma$ and $H_k>0$, $3\leq k\leq n$, the operators $L_j$ are elliptic for 
all $1\leq j\leq k-1$ and so $\mathcal{L}$ is elliptic as well. Furthermore, we observe that
\[
\mathrm{Tr}\mathcal{P}=\sum_{j=0}^{k-1}(1+\norm{\nabla u}^2)^{-\frac{k-j-1}{2}}\Big(\frac{h'}{h}(\inf_{\Sigma}u)\Big)^{k-j-1}\frac{c_{k-1}}{c_j}H_j. 
\]
Since $1/\sqrt{1+\norm{\nabla u}^2}\leq 1$, $h'/h(\inf_{\Sigma}u)<+\infty$ and, by the Newton inequalities
\[
 H_j\leq H_1^j<+\infty
\]
we conclude that $\mathcal{P}$ has trace bounded above and we can apply the Omori-Yau maximum principle for the operator $\mathcal{L}$. 
Hence, we can find a sequence $\{q_i\}_{i\in\enne}\subset\Sigma$ 
such that
\begin{equation}\label{OYLbelow}
(i)\ u(q_i)<u_*+\frac{1}{i},\ (ii) \ \norm{\nabla u(q_i)}<\frac{1}{i}, \ (iii) \ \mathcal{L}u(q_i)> -\frac{1}{i}.
\end{equation}
A straightforward computation using Proposition \ref{proplkuless} 
shows that
\begin{align*}
\mathcal{L}u\leq& -\frac{h'}{h}(u)\sum_{j=0}^{k-1}(1+\norm{\nabla u}^2)^{-\frac {k-1-j}{2}}\Big(\frac{h'}{h}(\inf_\Sigma u)\Big)^{k-j-1}\frac{c_{k-1}}{c_j}\pair{P_j\nabla u,\nabla u}\\
&-c_{k-1}\frac{h'}{h}(u)\frac{1}{(1+\norm{\nabla u}^2)^{(k-1)/2}}\Big(\frac{h'}{h}(\inf_\Sigma u)\Big)^{k-1}+\sqrt{1+\norm{\nabla u}^2}c_{k-1}H_k\\
&+c_{k-1}\sum_{j=1}^{k-1}(1+\norm{\nabla u}^2)^{-\frac {k-1-j}{2}}\Big(\frac{h'}{h}(\inf_\Sigma u)\Big)^{k-j-1}\Big(\frac{h'}{h}(\inf_{\Sigma}u)-\frac{h'}{h}(u)\Big)H_j.
\end{align*}
Evaluating the previous expression at $q_i$, using condition $(iii)$ in \eqref{OYLbelow} and taking the limit for $i\ra+\infty$, we 
find
\begin{displaymath}
0\leq -c_{k-1}\Big(\frac{h'}{h}(\inf_{\Sigma}u)\Big)^k+c_{k-1}\sup_{\Sigma}H_k
\end{displaymath}
and this concludes the proof. 
\end{proof}
\section{A Bernstein-type Theorem}
Recall now the Gauss equation
\begin{displaymath}
\R(X,Y)Z=(\overline{\R}(X,Y)Z)^T-\pair{AX,Z}AY+\pair{AY,Z}AX,
\end{displaymath}
for all tangent vector field $X,\ Y,\ Z\in T\Sigma$, where $(\overline{\R}(X,Y)Z)^T$ denotes the tangential component of $\overline{\R}(X,Y)Z$.  Hence, if $\{X,Y\}$ is any orthonormal basis of a tangent plane $\Pi\leq T_q\Sigma$, $q\in\Sigma$, the sectional curvature of $\Sigma$ is given by
\begin{align*}
K(X,Y)=&\overline{K}(X,Y)-\pair{AX,X}\pair{AY,Y}+\pair{AX,Y}^2\\
\geq & \overline{K}(X,Y)-\pair{AX,X}\pair{AY,Y}\\
\geq& \overline{K}(X,Y)-n^2H_1^2,
\end{align*}
where the last inequality follows by applying the Cauchy-Schwartz inequality.
In particular, if $M^{n+1}$ is a Lorentzian space form of constant sectional curvature $c$, then 
\begin{displaymath}
K(X,Y)\geq c-n^2H_1^2.
\end{displaymath}
Hence, if $\sup_{\Sigma}H_1<+\infty$ the Omori-Yau maximum principle holds on $\Sigma$ for semi-elliptic operators of the form $L=\text{Tr}(P\circ\Hess)$, where $P$ is a symmetric operator with trace bounded above. Applying the curvature estimates found in the previous section we are able to obtain the main result of this section, that extends Corollary 4.6 in \cite{aliashurtadopalmer} to spacelike hypersurfaces of constant higher order mean curvature. Notice that the previous estimates extend the ones given in \cite{aledoalias1} and \cite{aledoalias2}. Indeed, in this case the function $h$ has the expression
\begin{displaymath}
h(t)=
\left\{\begin{array}{ll}
\frac{1}{\sqrt{c}}\sinh(\sqrt{c}t) & \text{if $c>0$ and $t>0$}\\
t & \text{if $c=0$ and $t>0$}\\
\frac{1}{\sqrt{-c}}\sin(\sqrt{-c}t) & \text{if $c<0$ and $0<t<\pi/\sqrt{-c}$}
\end{array}
\right.
\end{displaymath}
Set $f_c(t)=h'(t)/h(t)$. Then 
\begin{displaymath}
f_c(t)=
\left\{\begin{array}{ll}
\sqrt{c}\coth(\sqrt{c}t) & \text{if $c>0$ and $t>0$}\\
\frac 1t & \text{if $c=0$ and $t>0$}\\
\sqrt{-c}\cot(\sqrt{-c}t) & \text{if $c<0$ and $0<t<\pi/\sqrt{-c}$}
\end{array}
\right.
\end{displaymath}
It is worth pointing out that $(f_c(t))^k$ is the $k$-mean curvature of the Lorentzian sphere of radius $t$ in the Lorentzian spaceform $M_c^{n+1}$ (when $\mathcal{I}^+(p)\neq\emptyset$), that is the level set
\[
\Sigma_c(t)=\{q\in\mathcal{I}^+(p)|d_p(q)=t\}.
\]
The following corollaries are straightforward.
\begin{corollary}\label{coroboundinfHk}
Let $M^{n+1}$ be an $n+1$- dimensional spacetime, $n \geq 3$, such that $K_{M}(\Pi)\leq c$, $c\in\erre$, for all timelike planes $\Pi$. Assume that there exists a point $p \in M$ such that $\mathcal{I}^+(p)\neq\emptyset$ and let $\psi:\Sigma^n\ra M^{n+1}$ be a spacelike hypersurface such that $\psi(\Sigma^n)\subset\mathcal{I}^+(p)\cap B^+(p,\delta)$ for some $\delta>0$ (with $\delta\leq\pi/\sqrt{-c}$ if $c<0$). Assume that either 
\begin{itemize}
\item[$(i)$] $k=2$ and $H_2$ is a positive function\\
or
\item[$(ii)$] $H_k$ is a positive function, $3\leq k\leq n$, and there exists an elliptic point $p_0\in\Sigma$.
\end{itemize}
Moreover, suppose that $\sup_{\Sigma}H_1<+\infty$ and that $\inf_\Sigma u<\pi/\sqrt{-c}$ if $c<0$. If the generalized Omori-Yau maximum principle holds on $\Sigma$, then
\begin{displaymath}
\inf_{\Sigma}H_k^{\frac 1k}\leq f_c\big(\sup_{\Sigma}u\big),
\end{displaymath}
where $u$ denotes the Lorentzian distance $d_p$ along the hypersurface.
\end{corollary}

\begin{corollary}\label{coroboundsupHk}
Let $M^{n+1}$ be an $n+1$- dimensional spacetime, $n \geq 3$, such that $K_{M}(\Pi)\geq c$, $c\in\erre$, for all timelike planes $\Pi$. Assume that there exists a point $p \in M$ such that $\mathcal{I}^+(p)\neq\emptyset$ and let $\psi:\Sigma^n\ra M^{n+1}$ be a spacelike hypersurface such that $\psi(\Sigma^n)\subset\mathcal{I}^+(p)$. Assume that either 
\begin{itemize}
\item[$(i)$] $k=2$ and $H_2$ is a positive function\\
or
\item[$(ii)$] $H_k$ is a positive function, $3\leq k\leq n$, and there exists an elliptic point $p_0\in\Sigma$.
\end{itemize}
Moreover, suppose that $\sup_{\Sigma}H_1<+\infty$ and that $\inf_\Sigma u<\pi/\sqrt{-c}$ if $c<0$. 
If the generalized Omori-Yau maximum principle holds on $\Sigma$, then
\begin{displaymath}
\sup_{\Sigma}H_k^{\frac 1k}\geq f_c\big(\inf_{\Sigma}u\big),
\end{displaymath}
where $u$ denotes the Lorentzian distance $d_p$ along the hypersurface.
\end{corollary}
Using the previous estimates we then obtain the following
\begin{theorem}
Let $M_c^{n+1}$ be a Lorentzian spaceform of constant sectional curvature $c$, $n\geq 3$, and let $p\in M_c^{n+1}$. Let $\Sigma$ be a complete spacelike hypersurface which is contained in $\mathcal{I}^+(p)$ such that either
\begin{itemize}
\item[$(i)$] $k=2$ and $H_2$ is a positive constant\\
or 
\item[$(ii$)] $H_k$ is constant, $3\leq k\leq n$, and there exists an elliptic point $p_0\in\Sigma$.
\end{itemize}
Moreover, assume that $\sup_{\Sigma}H_1<+\infty$. If $\Sigma$ is bounded from above by a level set of the Lorentzian distance function $d_p$ (with $d_p<\pi/\sqrt{-c}$ if $c<0$), then $\Sigma$ is necessarily a level set of $d_p$.
\end{theorem}
\begin{proof}
Our hypotheses imply that $\Sigma$ is contained in $\mathcal{I}^+(p)\cap B^+(p,\delta)$, with $\delta\leq\pi/\sqrt{-c}$ when $c<0$ 
and that $\Sigma$ has sectional curvature bounded from below. In particular the generalized Omori-Yau maximum principle holds on 
$\Sigma$ and we can apply Corollaries \ref{coroboundinfHk} and \ref{coroboundsupHk} to obtain
\begin{displaymath}
f_c(\sup_{\Sigma}u)\geq H_k^{\frac 1k}\geq f_c(\inf_{\Sigma}u).
\end{displaymath}
Hence, since $f_c$ is a decreasing function, $\sup_{\Sigma}u=\inf_{\Sigma}u=f_c^{-1}(H_k^{\frac 1k})$ and $\Sigma$ is necessarily the level set $d_p=f_c^{-1}(H_k^{\frac 1k})$.
\end{proof}

\begin{akn}
{The author would like to thank the referee for a 
careful reading of the original manuscript and for the useful suggestions.
}
\end{akn}

\bigskip

\bibliographystyle{amsplain}
\bibliography{biblio_comparison}

\end{document}